\documentstyle[amssymb,12pt]{amsart}

\newtheorem{lemma}{Lemma}
\newtheorem{thm}[lemma]{Theorem}

\newtheorem{cor}[lemma]{Corollary}

\newcommand{\DOT}{\setlength{\unitlength}{1pt}\begin{picture}(2.5,2)(1,1)
\put(1,2){\circle*{2}}\end{picture}}

\newcommand{\Calfdot}{{\cal{F}\!_{\DOT}\,}}

\newcommand{\Edot}{{E_{\DOT}\,}}
\newcommand{\Epdot}{{E_{\DOT}\!'\,}}

\newcommand{\Fdot}{{F\!_{\DOT}\,}}
\newcommand{\Fpdot}{{{F\!_{\DOT}}'\,}}

\newcommand{\Gdot}{{G_{\DOT}}}
\newcommand{\Gpdot}{{G_{\DOT}}\!'}
\newcommand{\Gppdot}{{G_{\DOT}}\!''}

\newcommand{\llra}{\relbar\joinrel\longrightarrow}
\newcommand{\lllra}{\relbar\joinrel\llra}

\newcommand{\rkm}{\stackrel{r[k,m]}{\relbar\joinrel\llra}}
\newcommand{\ckm}{\stackrel{c[k,m]}{\relbar\joinrel\llra}}
\newcommand{\cVIVI}{\stackrel{c[4,4]}{\relbar\joinrel\llra}}

\newcommand{\cnkm}{\stackrel{c[n{-}k,m]}{\relbar\joinrel\relbar\joinrel\lllra}}

\newcommand{\QED}{
\setlength{\unitlength}{1.0pt}%
\begin{picture}(7.5,7.5)
\put(0,-2.5){\rule{2.5pt}{2.5pt}}
\put(0,0){\rule{5pt}{2.5pt}}
\put(0,2.5){\rule{7.5pt}{2.5pt}}
\end{picture}\vspace{10pt}}

\newcommand{\Span}[1]{\langle #1 \rangle}

\begin{document}

\title{Pieri's Rule for Flag Manifolds and
 Schubert Polynomials}

\normalsize

\author{Frank Sottile}

\address{
        Department of Mathematics\\
        University of Toronto\\
        100 St.~George Street\\
	Toronto, Ontario  M5S 1A1\\
	Canada\\
	(416) 978-4031\\
	sottile\@math.toronto.edu
}
\date{2 May, 1995}

\maketitle
\normalsize
\begin{abstract}
We establish the formula for multiplication  by the class of
a special Schubert variety
in the integral cohomology ring of the flag manifold.
This formula also describes the multiplication of a Schubert
polynomial by either an elementary symmetric polynomial or a complete
homogeneous  symmetric polynomial.
Thus, we  generalize the classical Pieri's rule for symmetric
polynomials/Grassmann varieties to Schubert polynomials/flag manifolds.
Our primary technique is an explicit geometric description of certain
intersections of Schubert varieties.
This method allows us to compute additional structure constants
for the cohomology ring, which we express in terms of paths in the
Bruhat order on the symmetric group.
\end{abstract}

\section{Introduction}

Schubert polynomials had their origins in the study
of the cohomology of flag manifolds  by
Bernstein-Gelfand-Gelfand~\cite{BGG} and Demazure~\cite{Demazure}.
They were later defined by  Lascoux and
Sch\"utzenberger~\cite{Lascoux_Schutzenberger_polynomes_schubert},
who developed a purely combinatorial theory.

For each permutation $w$ in the symmetric group $S_n$ there is a Schubert
polynomial $\frak{S}_w$ in the variables $x_1,\ldots,x_{n-1}$.
When evaluated at certain Chern classes, a Schubert polynomial
gives the cohomology
class of a Schubert subvariety of the manifold of complete flags
in $\Bbb{C}\,^n$.
In this way, the
collection $\{\frak{S}_w\,|\, w\in S_n\}$ of Schubert polynomials
determines an integral basis for the cohomology ring of
the flag manifold.
Thus there exist integer structure constants $c^u_{w\,v}$ such that
$$
\frak{S}_w\cdot\frak{S}_v = \sum_u  c^u_{w\,v}\frak{S}_u.
$$
No formula is known, or even conjectured, for these constants.
There are, however, a few special cases in which they are known.

One important case is Monk's rule~\cite{Monk}, which characterizes
the algebra of Schubert polynomials.
While this is usually attributed to Monk,
Chevalley simultaneously established the analogous formula
for generalized flag manifolds in a manuscript that was only
recently published~\cite{Chevalley91}.
Let $t_{k\,k+1}$ be the transposition interchanging $k$ and $k+1$.
Then $\frak{S}_{t_{k\,k+1}} =  x_1{+}\cdots{+}x_k=s(x_1,\ldots,x_k)$,
the first elementary symmetric polynomial.
For any permutation $w\in S_n$, Monk's rule states
$$
\frak{S}_w \cdot \frak{S}_{t_{k\,k+1}} \
=\ \frak{S}_w \cdot s_1(x_1,\ldots,x_k)
\ =\ \sum \frak{S}_{w t_{a\,b}},
$$
where $t_{a\,b}$ is the transposition interchanging $a$ and $b$, and
the sum is over all $a\leq k<b$ where $w(a)<w(b)$ and if
$a<c<b$, then $w(c)$ is not between $w(a)$ and $w(b)$.
\smallskip

The classical Pieri's rule computes the product of
a Schur polynomial by either a complete homogeneous
symmetric polynomial or an elementary symmetric polynomial.
Our main result is a formula for Schubert polynomials
and the cohomology of flag manifolds which generalizes both Monk's rule
and the classical Pieri's rule.

Let $s_m(x_1,\ldots,x_k)$ and $s_{1^m}(x_1,\ldots,x_k)$ be
respectively the complete homogeneous
and elementary symmetric polynomials of degree
$m$ in the variables $x_1,\ldots,x_k$ and let $\ell(w)$ be the length of a
permutation $w$.
These polynomials are the cohomology classes of special Schubert varieties.
We will show
\medskip

\noindent{\bf Theorem~\ref{thm:main}.}
{\em Let $k,m,n$ be positive integers, and let $w\in S_n$.
\begin{enumerate}
\item[I.]
$\frak{S}_w\cdot  s_m(x_1,\ldots,x_k) =
\sum_{w'} \frak{S}_{w'}$,
the sum over all $w' = w t_{a_1\,b_1}\cdots t_{a_m\,b_m}$, where
$a_i\leq k < b_i$  and
$\ell(w t_{a_1\,b_1}\cdots t_{a_i\,b_i}) =  \ell(w) + i$
for $1\leq i\leq m$ with the integers
$b_1,\ldots, b_m$ distinct.
\item[II.]
$\frak{S}_w\cdot  s_{1^m}(x_1,\ldots,x_k) =
\sum_{w'} \frak{S}_{w'}$,
the sum over all $w'$ as in {\em I}, except that now
the integers $a_1,\ldots,a_m$ are distinct.
\end{enumerate}
}
\medskip

Both  $s_m(x_1,\ldots,x_k)$ and $s_{1^m}(x_1,\ldots,x_k)$
are Schubert polynomials, so Theorem~\ref{thm:main} computes some
of the structure constants in the cohomology ring of the flag manifold.

These formulas were stated in a different form by Lascoux and Sch\"utzenberger
in~\cite{Lascoux_Schutzenberger_polynomes_schubert}, where an
algebraic proof was suggested.
They were later independently conjectured in yet another form by Bergeron and
Billey~\cite{Bergeron_Billey}.
Our formulation facilitates our proofs.
Using geometry, we expose a surprising connection to the classical
Pieri's rule, from which we deduce Theorem~\ref{thm:main}.
These methods enable the determination of additional structure constants.
We further generalize Theorem~\ref{thm:main} to give a formula for the
multiplication of a Schubert polynomial  by a  hook Schur polynomial,
indicating a relation between multiplication of Schubert
polynomials and paths in the Bruhat order in $S_n$.

This exposition is organized as follows:
Section 2 contains preliminaries about Schubert polynomials
while Section 3 is devoted to the flag manifold.
In Section 4 we deduce our main results from
a geometric lemma proven in Section~5.
Two examples are described in Section~6, illustrating the
geometry underlying the results of Section~5.
We remark that while our results are stated in terms of the
integral cohomology of the complex manifold of complete flags,
our results and proofs are valid for the Chow rings of
flag varieties defined over any field.

We would like to thank Nantel Bergeron and
Sara Billey for suggesting these problems
and Jean-Yves Thibon for showing us the work of Lascoux
and Sch\"utzenberger.

\section{Schubert Polynomials}

In~\cite{BGG,Demazure} cohomology classes of Schubert subvarieties of
the flag manifold were obtained from the class of a point using
repeated correspondences in $\Bbb{P}^1$-bundles, which may be described
algebraically as ``divided differences.''
Subsequently, Lascoux and
Sch\"utzenberger~\cite{Lascoux_Schutzenberger_polynomes_schubert}
found explicit polynomial representatives
for these classes.
We outline Lascoux and Sch\"utzenberger's construction of Schubert
polynomials.
For a more complete account see~\cite{Macdonald_schubert}.

For an integer $n>0$, let $S_n$ be the group of
permutations of $[n] = \{1,2,\ldots,n\}$.
Let $t_{a\, b}$ be the transposition interchanging $a <b$.
Adjacent transpositions $s_i = t_{i\,i{+}1}$
generate $S_n$.
The {\em length} $\ell(w)$ of a
permutation $w$ is the minimal length of a factorization
into adjacent transpositions.
If $w = s_{a_1} s_{a_2}\cdots s_{a_m}$ is such a factorization, then
the sequence $(a_1,\ldots,a_m)$ is a {\em reduced word} for $w$.
The length of $w$ also counts the  inversions of $w$,
those pairs $i<j$ where $w(i)>w(j)$.
It follows that  $\ell(wt_{a\, b}) = \ell(w){+}1$ if and only if
$w(a)<w(b)$ and whenever  $a<c<b$, either $w(c)< w(a)$ or
$w(b)<w(c)$.

For each integer $n>1$, let $R_n = \Bbb{Z}[x_1,\ldots,x_n]$.
The group $S_n$ acts on $R_n$  by permuting the variables.
Let $f\in R_n$  and  let $s_i$ be an adjacent transposition.
The polynomial
$f - s_i f$ is antisymmetric in $x_i$ and $x_{i+1}$, and so is
divisible by $x_i - x_{i+1}$.
Thus we may define the linear  divided difference operator
$$
\partial_i = (x_i-x_{i+1})^{-1} (1 - s_i).
$$
If $f$ is symmetric in $x_i$ and $x_{i+1}$, then $\partial_i f$ is zero.
Otherwise $\partial_i f$ is symmetric in $x_i$ and $x_{i+1}$.
Divided differences satisfy
\begin{eqnarray*}
\partial_i\circ \partial_i & = & 0 \\
\partial_i \circ\partial_j &=& \partial_j \circ \partial_i \ \
\ \ \ \ \ \ \ \ \ \ \mbox{ if } |i-j|\geq 2\\
\qquad\qquad\partial_{i+1}\circ\partial_i\circ\partial_{i+1} & = &
\partial_i\circ\partial_{i+1}\circ\partial_i
\end{eqnarray*}
It follows that if $(a_1,\ldots,a_p)$ is a reduced word for a permutation
$w$, the composition of divided differences
$\partial_{a_1}\circ\cdots\circ\partial_{a_p}$
depends only upon $w$ and not upon the reduced word chosen.
This defines an operator $\partial_{w}$ for each $w\in S_n$.

Let $w_0$ be the longest permutation in $S_n$, that is
$w_0(j) = n{+}1{-}j$.
For $w \in S_n$, define the {\em Schubert polynomial} $\frak{S}_{w}$
by
$$
\frak{S}_{w} = \partial_{w^{-1}w_0}
\left( x_1^{n-1} x_2^{n-2}\cdots x_{n-1} \right).
$$
The degree of $\partial_i$ is $-1$, so $\frak{S}_{w}$ is homogeneous of
degree $ {n\choose 2} - \ell(w^{-1}w_0) = \ell(w)$.

Let $\cal{S}\subset R_n$ be the ideal generated by the
non-constant symmetric polynomials.
The set  $\{\frak{S}_{w}\,|\, w\in S_n\}$ of Schubert polynomials
is a basis for
$\Bbb{Z}\{ x_1^{i_1}\cdots x_{n-1}^{i_{n-1}}\,|\, i_j \leq n{-}j \}$,
a transversal to $\cal{S}$ in $R_n$.
Thus Schubert polynomials are  explicit polynomial
representatives of an integral basis for the ring $H_n = R_n/\cal{S}$.
Courting ambiguity, we will use the same notation for
Schubert polynomials in $R_n$ as for their images in the rings
$H_n$.

Recently, other descriptions have been discovered
for Schubert
polynomials~\cite{Bergeron,BJS,Fomin_Kirillov,Fomin_Stanley}.
Combinatorists often define Schubert polynomials $\frak{S}_w$ for all
$w\in S_\infty = \cup_{n=1}^\infty S_n$.
One may show that our results are valid in this wider context.
\smallskip

A {\em partition} $\lambda$ is a decreasing sequence
$\lambda_1 \geq \lambda_2\geq\cdots\geq \lambda_k$
of positive integers, called the {\em parts} of $\lambda$.
Given a partition $\lambda$ with at most $k$ parts,
one may define a Schur polynomial
$s_\lambda = s_\lambda(x_1,\ldots,x_k)$, which is a symmetric
polynomial in the variables $x_1,\ldots,x_k$.
For a more complete treatment of symmetric polynomials and
Schur polynomials, see~\cite{Macdonald_symmetric}.

The collection of Schur polynomials  forms an
integral basis for the ring of symmetric polynomials,
 $\Bbb{Z}[x_1,\ldots,x_k]^{S_k}$.
The Littlewood-Richardson rule is a formula for the
structure constants $c^\lambda_{\mu\nu}$ of this ring, called
{\em Littlewood-Richardson coefficients} and defined by
$$
s_\mu \cdot s_\nu \ =\ \sum_\lambda \, c^\lambda_{\mu\nu}\, s_\lambda.
$$
If $\lambda$ and $\mu$ are partitions satisfying $\lambda_i \geq \mu_i$
for all $i$, we write $\lambda \supset \mu$.
This defines a partial order on the collection of partitions,
called Young's lattice.
Since $c^\lambda_{\mu\nu} = 0$ unless $\lambda \supset \mu$ and
$\lambda\supset \nu$~(cf. \cite{Macdonald_symmetric}), we see that
$\cal{I}_{n,k}=\{ s_\lambda\,|\, \lambda_1 \geq n-k\}$ is an ideal.
Let $A_{n,k}$ be the quotient ring
$\Bbb{Z}[x_1,\ldots,x_k]^{S_k}/\cal{I}_{n,k}$.

To a partition $\lambda$ we may associate its Young diagram,
also denoted $\lambda$, which is a left-justified
array of boxes in the plane
with $\lambda_i$ boxes in the $i$th  row.
If $\lambda \supset \mu$, then the Young diagram of $\mu$ is
a subset of that of $\lambda$, and the skew diagram $\lambda/\mu$
is the set theoretic difference $\lambda-\mu$.
If each column of $\lambda/\mu$ is either empty or a single box, then
$\lambda/\mu$ is a {\em skew row} of {\em length} $m$, where $m$
is the number of boxes in $\lambda/\mu$.
The transpose $\mu^t$ of a partition $\mu$ is the partition whose Young
diagram is the transpose of that of $\mu$.
We call the  transpose of  a skew row a {\em skew column}.
The map defined by $s_\lambda \mapsto s_{\lambda^t}$ is a ring
isomorphism $A_{n,k} \rightarrow A_{n,n-k}$.

For example, let $\lambda = (5,2,1)$ and $\mu = (3,1)$
then $\lambda/\mu$ is a skew row of length 4 and
$\mu^t = (2,1,1)$.
The following are the Young diagrams of $\lambda$, $\mu$,
$\lambda/\mu$, and $\mu^t$:

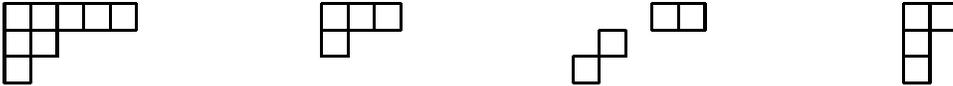
\begin{picture}(400,60)
\put(20,10){\begin{picture}(50,30)
\thicklines
\put(0, 0){\line(1,0){10}}
\put(0,10){\line(1,0){20}}
\put(0,20){\line(1,0){50}}
\put(0,30){\line(1,0){50}}
\put( 0, 0){\line(0,1){30}}
\put(10, 0){\line(0,1){30}}
\put(20,10){\line(0,1){20}}
\put(30,20){\line(0,1){10}}
\put(40,20){\line(0,1){10}}
\put(50,20){\line(0,1){10}}

\end{picture}}

\put(140,10){\begin{picture}(50,30)
\thicklines
\put( 0,10){\line(1,0){10}}
\put( 0,20){\line(1,0){30}}
\put( 0,30){\line(1,0){30}}
\put( 0,10){\line(0,1){20}}
\put(10,10){\line(0,1){20}}
\put(20,20){\line(0,1){10}}
\put(30,20){\line(0,1){10}}
\end{picture}}

\put(235,10){\begin{picture}(50,30)
\thicklines
\put( 0, 0){\line(1,0){10}}
\put( 0,10){\line(1,0){20}}
\put(10,20){\line(1,0){10}}
\put( 0, 0){\line(0,1){10}}
\put(10, 0){\line(0,1){20}}
\put(20,10){\line(0,1){10}}

\put(30,20){\line(1,0){20}}
\put(30,30){\line(1,0){20}}
\put(30,20){\line(0,1){10}}
\put(40,20){\line(0,1){10}}
\put(50,20){\line(0,1){10}}
\end{picture}}

\put(360,10){\begin{picture}(40,30)
\thicklines
\put(0, 0){\line(1,0){10}}
\put(0,10){\line(1,0){10}}
\put(0,20){\line(1,0){20}}
\put(0,30){\line(1,0){20}}
\put( 0, 0){\line(0,1){30}}
\put(10, 0){\line(0,1){30}}
\put(20,20){\line(0,1){10}}
\end{picture}}

\end{picture}

If $w$ has only one {\em descent} ($k$ such that $w(k) >
w(k{+}1)$), then $w$ is said to be {\em Grassmannian}
of descent $k$ and $\frak{S}_{w}$ is the Schur polynomial
$s_{\lambda}(x_1,\ldots,x_k)$.
Here $\lambda$ is the {\em shape} of $w$, the partition
with $k$ parts where $\lambda_{k+1-j} = w(j){-}j$.
For integers $k,m$ define
$r[k,m]$ and $c[k,m]$ to be the Grassmannian permutations of descent $k$
with shapes  $(m,0,\ldots,0) = m$ and
$(1^m,0,\ldots,0) = 1^m$, respectively.
These are the $m+1$-cycles
\begin{eqnarray*}
r[k,m] &=&
(k{+}m\,\,\,\,k{+}m{-}1\,\ldots\,k{+}2\,\,\,\,k{+}1\,\,\,\,k)\\
c[k,m]  &=&
(k{-}m{+}1\,\,\,\,k{-}m{+}2\,\ldots\,k{-}1\,\,\,\,k\,\,\,\,k{+}1).
\end{eqnarray*}

\section{The Flag Manifold}

Let $V$ be an $n$-dimensional complex vector space.
A {\em flag} $\Fdot$ in $V$ is a sequence
$$
\{0\}\  =\ F_0 \subset F_1 \subset F_2\subset \cdots \subset F_{n-1}
\subset F_n\ =\ V,
$$
of linear subspaces with $\dim_{\Bbb{C}} F_i = i$.
The set of all flags is a $\frac{1}{2}n(n-1)$ dimensional complex
manifold, called the flag manifold
 and denoted $\Bbb{F}(V)$.
Over $\Bbb{F}(V)$, there is a tautological flag $\Calfdot$ of
bundles whose fibre at a point $\Fdot$ is the flag $\Fdot$.
Let $x_i$  be the Chern class of the line bundle
$\cal{F}_i/\cal{F}_{i-1}$.
Then the integral cohomology ring of $\Bbb{F}(V)$ is
$H_n = \Bbb{Z}[x_1,\ldots,x_n]/\cal{S}$, where $\cal{S}$
is the ideal generated by those non-constant polynomials
which are symmetric in
$x_1,\ldots,x_n$.
This description is due to Borel~\cite{Borel}.

Given a subset $S \subset V$, let $\Span{S}$ be its linear span
and for linear subspaces $W\subset U$ let $U-W$ be their
set theoretic difference.
An ordered
basis $f_1,f_2,\ldots,f_n$ for $V$ determines a flag $\Edot$; set
$E_i = \Span{f_1,\ldots,f_i}$ for $1\leq i \leq n$.
In this case, write $\Edot = \Span{f_1,\ldots,f_n}$
and call $f_1,\ldots,f_n$ a {\em basis} for $\Edot$.
A fixed flag $\Fdot$ gives a decomposition due to
Ehresmann~\cite{Ehresmann} of $\Bbb{F}(V)$ into
affine cells indexed by permutations $w$ of $S_n$.
The cell determined by $w$ is
$$
X^{\circ}_w \Fdot \  = \  \{ \Edot=\Span{f_1,\ldots,f_n}\,|\,
f_i \in F_{n+1-w(i)}-F_{n-w(i)}, \,1\leq i\leq n\}.
$$
The complex codimension of $X^{\circ}_w \Fdot$  is $\ell(w)$
and its closure is the Schubert subvariety $X_w\Fdot$.
Thus the cohomology ring of $\Bbb{F}(V)$
has an integral basis given by the
cohomology classes\footnote{Strictly speaking, we mean the
classes Poincar\'e dual to the fundamental cycles in homology.}
$[X_w\Fdot]$ of the Schubert subvarieties.
That is, $H^*\Bbb{F}(V) = \bigoplus_{w\in S_n}{\Bbb Z} [X_w\Fdot]$.

Independently,  Bernstein-Gelfand-Gelfand~\cite{BGG}
and Demazure~\cite{Demazure}
related this description to Borel's, showing
$[X_w\Fdot] = \partial_{w^{-1}w_0}[\{\Fdot\}]$.
Later, Lascoux and
Sch\"utzenberger~\cite{Lascoux_Schutzenberger_polynomes_schubert}
obtained  polynomial
representatives $\frak{S}_w$ for $[X_w\Fdot]$
by choosing
$x_1^{n-1} x_2^{n-2}\cdots x_{n-1}$ for the
representative of the class $[\{\Fdot\}] = \frak{S}_{w_0}$ of a point.
We use the term Schubert polynomial for both the polynomial and
the associated cohomology class.

This Schubert polynomial basis for cohomology diagonalizes
the intersection pairing; If
$\ell(w) + \ell(w') = \dim\Bbb{F}(V) = \frac{1}{2}n(n-1)$, then
$$
\frak{S}_w\cdot \frak{S}_{w'} =
\left\{
\begin{array}{ll} \frak{S}_{w_0} & \mbox{ if } w' = w_0 w\\
0 & \mbox{ otherwise}
\end{array} \right.
$$
\smallskip

For each $k\leq \dim V =n$, the set of all $k$-dimensional subspaces of $V$
is a $k(n{-}k)$ dimensional complex manifold, called the Grassmannian of
$k$-planes in $V$, written  $G_kV$.
The cohomology ring of $G_kV$ is a quotient of the ring of symmetric
polynomials in the Chern roots $x_1,\ldots,x_k$ of its tautological
$k$-plane bundle.
This identifies it with the ring $A_{n,k}$ of Section 2.

A fixed flag $\Fdot$ gives a decomposition of $G_kV$ into
cells indexed by partitions $\lambda$ with $k$ parts, none
exceeding $n{-}k$.
The closure of such a cell is the Schubert variety
$$
\Omega_\lambda \Fdot =
\{ H \in G_kV \,|\, \dim H\cap F_{n-k+j-\lambda_j} \geq j
\mbox{ for } 1\leq j\leq k\},
$$
whose codimension is
$\lambda_1{+}{\cdots}{+}\lambda_k = |\lambda|$.
The classes $[\Omega_\lambda\Fdot]$  form a basis for the
cohomology ring of $G_kV$ and $[\Omega_\lambda\Fdot]$
is the Schur polynomial $s_\lambda(x_1,\ldots,x_k)$.
We use the term Schur polynomial for both the polynomial and its
image in the cohomology ring of $G_kV$.

The Schur polynomial $s_{m}$ is  the complete homogeneous
symmetric polynomial of degree $m$ in $x_1,\ldots,x_k$.
The Schur polynomial $s_{1^m}$ is  the $m$th
elementary symmetric polynomial in $x_1,\ldots,x_k$.
Pieri's rule is a formula for multiplying Schur polynomials by
either $s_m$ or $s_{1^m}$.
For  $s_m$, it states
$$
s_\mu \cdot s_m\  = \
\sum s_{\lambda},
$$
the sum over all partitions $\lambda$ with
$n{-}k\geq\lambda_1\geq\mu_1\geq\cdots\geq\lambda_k\geq\mu_k$ and
$|\lambda| =  m{+} |\mu|$.
That is, those partitions $\lambda\supset \mu$ with $\lambda/\mu$ a skew row
of length $m$.

To obtain the analogous formula for $s_{1^m}$, use the isomorphism
$A_{n,k} \rightarrow A_{n,n-k}$ given by
$s_\lambda \mapsto s_{\lambda^t}$.
Doing so, we see that
$$
s_\mu  \cdot s_{1^m}\  = \ \sum s_\lambda,
$$
the sum over all partitions $\lambda$ with
$\lambda\supset \mu$ with $(\lambda/\mu)^t$ is a skew row of length $m$.
That is, those $\lambda\supset \mu$ with $\lambda/\mu$ a skew column
of length $m$.
\smallskip

If $Y\subset V$ has codimension $d$, then $G_kY \subset G_kV$ is a
Schubert subvariety whose indexing partition is $d^k$, the partition
with $k$ parts each equal to $d$.
It follows that $\Omega_{(n{-}k)^k}\Fdot = \{F_k\}$, so $s_{(n{-}k)^k}$
is the class of a point.

The basis of Schur polynomials diagonalizes the intersection pairing;
For a partition $\lambda$, let $\lambda^c$ be the partition
$(n{-}k{-}\lambda_k,{\ldots},n{-}k{-}\lambda_1)$.
If $|\mu| {+}|\lambda| = k(n{-}k)$, then
$$
s_\lambda \cdot s_\mu =
\left\{
\begin{array}{ll} s_{(n{-}k)^k}& \mbox{ if } \lambda^c = \mu\\
0 & \mbox{ otherwise }
\end{array}
\right ..
$$
We use this to  reformulate Pieri's rule.
Suppose $|\mu|+|\lambda|+m = k(n-k)$, then
$$
s_\mu \cdot s_{\lambda^c}\cdot s_m = \left\{
\begin{array}{ll} s_{(n-k)^k} &
\mbox{ if $\lambda/\mu$ is a skew row of length $m$}\\
0 &\mbox{ otherwise}\end{array}\right. .
$$
\smallskip

For $k\leq n$, the association $\Edot \mapsto E_k$
defines a map $\pi :\Bbb{F}(V) \rightarrow G_kV$.
The functorial map $\pi^*$ on
cohomology is simply the inclusion into $H_n$
of polynomials symmetric in $x_1,\ldots,x_k$.
That is, $A_{n,k} \hookrightarrow H_n$.
If $\lambda$ is a partition with $k$ parts and $w$
the Grassmannian permutation of descent $k$ and shape $\lambda$,
then $\pi^* s_\lambda = \frak{S}_w$.

Under the Poincar\'e duality isomorphism between homology and
cohomology groups, the functorial map $\pi_*$ on homology induces a
a group homomorphism $\pi_*$ on cohomology.
While $\pi_*$ is not a ring homomorphism, is does satisfy
the projection formula (see Example 8.17 of~\cite{Fulton_intersection}):
$$
\pi_*(\alpha\cdot \pi^* \beta) = (\pi_* \alpha)\cdot \beta,
$$
where $\alpha$ is a cohomology class on $\Bbb{F}(V)$ and $\beta$
is a cohomology class on $G_kV$.

\section{Pieri's Rule for Flag Manifolds}

An open problem  is to find  the analog
of the Littlewood-Richardson rule for Schubert polynomials.
That is, determine the structure constants $c^u_{w\,v}$ for
the Schubert basis of the cohomology of flag manifolds, which are
defined by
\begin{equation} \label{eq:structure}
\frak{S}_w \cdot \frak{S}_v =
\sum_u c^u_{w\,v} \frak{S}_u.
\end{equation}
These constants are positive integers as they count
the points in a suitable triple intersection
of Schubert subvarieties.
They are are known only in some special cases.

For example, if both $w$ and $v$ are Grassmannian
permutations of descent  $k$ so that $\frak{S}_w$ and
$\frak{S}_v$ are
symmetric polynomials in the variables $x_1,\ldots,x_k$, then
(\ref{eq:structure}) is
the classical Littlewood-Richardson rule.

Another  case is Monk's rule, which states:
$$
	\frak{S}_w\cdot \frak{S}_{t_{k\,k{+}1}}
	= \sum \frak{S}_{w t_{a\,b}},
$$
the sum over all $a\leq k <b$ with
$\ell(w t_{a\,b})=\ell(w)+1$.
The Schubert polynomial $\frak{S}_{t_{k\,k{+}1}}$ is
$s_1(x_1,\ldots,x_k)$.
We use geometry to generalize this formula, giving an analog of
the classical Pieri's rule.
\smallskip

Let $w,w' \in S_n$.
Write $w \rkm w'$ if
there exist integers $a_1,b_1,\ldots,a_m,b_m$ with
\begin{enumerate}
\item $a_i\leq k <b_i$ for $1\leq i\leq m$ and
$w' = wt_{a_1\,b_1}\cdots t_{a_m\,b_m}$,
\item $\ell(w t_{a_1\,b_1}\cdots t_{a_i\,b_i}) = \ell(w) +i$, and
\item the integers $b_1, b_2,\ldots, b_m$ are distinct.
\end{enumerate}
Similarly, $w\ckm w'$
if we have integers $a_1,\ldots,b_m$ as in (1) and (2) where now
\begin{enumerate}
\item[(3)$'$] the integers $a_1,a_2,\ldots, a_m$ are distinct.
\end{enumerate}
Our primary result is the following.

\begin{thm} \label{thm:main}
Let $w \in S_n$.
Then
\begin{enumerate}

\item[I.] For all $k$ and $m$ with $k+m \leq n$,  we have \ \
${\displaystyle
\frak{S}_{w}\cdot \frak{S}_{r[k,m]}
= \sum_{w \rkm w'}
\frak{S}_{w'}}$.
\smallskip

\item[II.] For all $ m\leq k\leq n$, we have\rule{0pt}{20pt} \ \
${\displaystyle
\frak{S}_{w}\cdot \frak{S}_{c[k,m]}
= \sum_{w \ckm w'}
\frak{S}_{w'}}$.

\end{enumerate}
\end{thm}

Theorem~\ref{thm:main} may be alternatively stated in terms of the
structure constants $c^u_{w\,v}$.
\medskip

\noindent{\bf Theorem 1$'\!$.} \ {\em Let $w, w' \in S_n$.
Then
\begin{enumerate}

\item[I.] For all integers $k,m$ with $k+m\leq n$, \ \
${\displaystyle
c^{w'}_{w\, r[k,m]} = \left\{\begin{array}{ll} 1 &\mbox{ if } w\rkm w'\\
0 & \mbox{ otherwise}\end{array}\right. 
}$.

\item[ II.] For all integers $k,m$ with $m\leq k\leq n$,\rule{0pt}{28pt} \ \
${\displaystyle
c^{w'}_{w\, c[k,m]} = \left\{\begin{array}{ll} 1 &\mbox{ if } w\ckm w'\\
0 & \mbox{ otherwise}\end{array}\right. 
}$.

\end{enumerate}
}
\bigskip

We first show the equivalence of parts I and II and then establish
part I.
An order $<_k$ on
$S_n$ is introduced, and we show that $c^{w'}_{w\, r[k,m]}$ is 0 unless
$w<_k w'$.
A geometric lemma enables us to compute
$c^{w'}_{w\, r[k,m]}$ when $w<_k w'$.

\begin{lemma}\label{lemma:equivalent}
Let  $w_0$ be the longest permutation in  $S_n$,
and $k{+}m \leq n$.
Then
\begin{enumerate}
\item $w_0 r[k,m] w_0 = c[n{-}k,m]$.
\item Let $w, w' \in S_n$.   Then
$w\rkm w'$ if and only if $\,\,\,w_0 w w_0\cnkm w_0w'w_0$.
\item The map induced by $\frak{S}_w \mapsto \frak{S}_{w_0ww_0}$
is an automorphism of $H_n$.
\item Statements {\em I} and {\em II}
of Theorem~\ref{thm:main}$'$ are equivalent.
\end{enumerate}
\end{lemma}

This automorphism $\frak{S}_w \mapsto \frak{S}_{w_0ww_0}$
is the Schubert polynomial analog of the map
$s_\lambda(x_1,\ldots,x_k) \mapsto
s_{\lambda^t}(x_1,\ldots,x_{n-k})$ for Schur polynomials.
\medskip

\noindent{\bf Proof:}
Statements (1) and (2) are easily verified, as $w_0(j) = n+1-j$.

Statement (3) is also immediate, as
$\frak{S}_w \mapsto \frak{S}_{w_0ww_0}$
leaves Monk's rule invariant and Monk's rule characterizes
the algebra of Schubert polynomials.

For (4), suppose $k+m \leq n$ and $w, w' \in S_n$ and
let $\overline{w}$ denote $w_0ww_0$.
The isomorphism $\frak{S}_v \mapsto \frak{S}_{\overline{v}}$
of (3) shows $c^{w'}_{w\, r[k,m]}=
c^{\overline{w'}}_{\overline{w}\,\overline{r[k,m]}}$.
Part (1) shows
$c^{\overline{w'}}_{\overline{w}\,\overline{r[k,m]}} =
c^{\overline{w'}}_{\overline{w}\,c[n{-}k,m]}$.
Then (2) shows the equality of the two statements of
Theorem~\ref{thm:main}$'$.
\QED

Let $<_k$ be the transitive closure of the relation given by
$w <_k w'$, whenever $w' = w t_{a\,b}$ with $a\leq k<b$ and
$\ell(w\, t_{a\,b}) = \ell(w){+}1$.
We call $<_k$ the {\em $k$-Bruhat order},
in~\cite{Lascoux_Schutzenberger_Symmetry}
it is the $k$-colored Ehresmano\"edre.

\begin{lemma}\label{lemma:order}
If $\,c^{w'}_{w\, r[k,m]} \neq 0$, then $w<_k w'$ and
$\ell(w') = \ell(w) + m$.
\end{lemma}

\noindent{\bf Proof:}
By Monk's rule, $w<_kw'$ if and only if $\frak{S}_{w'}$ appears with a
non-zero coefficient when
$\frak{S}_w (\frak{S}_{t_{k\,k{+}1}})^{\ell(w')-\ell(w)}$
is written as a sum of Schubert polynomials.

Since $r[k,m] = t_{k\,k{+}1} \cdot t_{k\,k{+}2}\cdots t_{k\,k{+}m}$,
Monk's rule  shows that $\frak{S}_{r[k,m]}$ is a summand of
$(\frak{S}_{t_{k\,k+1}})^m$ with coefficient 1.
Thus the coefficient of $\frak{S}_{w'}$ in the expansion of
$\frak{S}_w \cdot (\frak{S}_{t_{k\,k+1}})^m$  exceeds the
coefficient of $\frak{S}_{w'}$ in $\frak{S}_w \cdot \frak{S}_{r[k,m]}$.
Hence $c_{w\,r[k,m]}^{w'} = 0$ unless
$w<_k w'$ and $\ell(w') = \ell(w) +m$.
\QED

In Section 5 we use geometry to prove the following lemma.

\begin{lemma}\label{lemma:pushforward}
 Let $w<_k w'$ be permutations in $S_n$.
Suppose $w' = w t_{a_1\,b_1}\cdots t_{a_m\,b_m}$,
where $a_i\leq k<b_i$, and
$\ell(w t_{a_1\,b_1}\cdots t_{a_i\,b_i}) = \ell(w)+i$.
Let $d = n-k-\#\{b_1,\ldots,b_m\}$.
Then
\begin{enumerate}
\item There is a cohomology class $\delta$ on $G_kV$
such that
$\pi_*(\frak{S}_w\cdot \frak{S}_{w_0w'}) = \delta \cdot s_{d^k}$.
\item If $w\rkm w'$, then there are partitions
$\lambda \supset \mu$ where $\lambda/\mu$ is a skew row of length $m$
whose $j$th row has  length $\#\{i\,|\, a_i = j\}$
and
$\pi_*(\frak{S}_w\cdot \frak{S}_{w_0w'}) = s_{\mu}\cdot s_{\lambda^c}$.
\end{enumerate}
\end{lemma}

\noindent{\bf Proof of Theorem~\ref{thm:main}$'$:}
By Lemma~\ref{lemma:order}, we need only show that if $w<_k w'$
and $\ell(w')-\ell(w) = m$, then
$$
c^{w'}_{w\,r[k,m]} =
\left\{ \begin{array}{ll} 1 &\mbox{ if } w \rkm w'\\
0 & \mbox{ otherwise} \end{array} \right. .
$$

Begin by multiplying the identity $\frak{S}_w\cdot \frak{S}_{r[k,m]} =
\sum_v\, c^v_{w\, r[k,m]}\, \frak{S}_v$
by $\frak{S}_{w_0\, w'}$ and use the intersection
pairing to obtain
$$
\frak{S}_w\cdot\frak{S}_{w_0\, w'}\cdot\frak{S}_{r[k,m]}
\ =\  c^{w'}_{w\,r[k,m]}\,  \frak{S}_{w_0}.
$$
Recall that $\frak{S}_{r[k,m]} = \pi^* s_m(x_1,\ldots,x_k)$.
As  $\frak{S}_{w_0}$ and $s_{(n-k)^k}$ are the classes of  points,
$\pi_*\frak{S}_{w_0} = s_{(n-k)^k}$.
Apply the map $\pi_*$  and then the
projection formula to obtain:
\begin{eqnarray*}
\pi_*(\frak{S}_w\cdot\frak{S}_{w_0\, w'}\cdot \pi^* s_m)
&=& c^{w'}_{w\,r[k,m]}\, \pi_*( \frak{S}_{w_0})\\
\pi_*(\frak{S}_w\cdot\frak{S}_{w_0\, w'}) \cdot s_m
&=& c^{w'}_{w\,r[k,m]} \, s_{(n-k)^k}.
\end{eqnarray*}
By part (1) of Lemma~\ref{lemma:pushforward}, there is a cohomology
class $\delta$ on $G_kV$ with
$$
\pi_*(\frak{S}_w\cdot\frak{S}_{w_0\, w'}) \cdot s_m\
= \ \delta \cdot  s_{d^k} \cdot s_m
$$
But $s_{d^k} \cdot s_m = 0$ unless $d+m \leq n-k$.
Since $d = n-k-\#\{b_1,\ldots,b_m\}\geq n-k-m$,
we see that $ c^{w'}_{w\,r[k,m]} =0$ unless
$m = \#\{b_1,\ldots,b_m\}$, which implies
$w \rkm w'$.

To complete the proof of Theorem~\ref{thm:main}$'$, suppose that
$w \rkm w'$.
By part (2) of Lemma~\ref{lemma:pushforward},
there are partitions
$\lambda \supset \mu$ with $\lambda/\mu$ a skew row of length $m$
where we have
$\pi_*(\frak{S}_w\cdot \frak{S}_{w_0w'}) = s_{\mu}\cdot s_{\lambda^c}$.
Then
$$
 \pi_*(\frak{S}_w\cdot\frak{S}_{w_0\, w'}) \cdot s_m\
= \ s_{\mu}\cdot s_{\lambda^c} \cdot s_m \
= \ s_{(n-k)^k},
$$
by the ordinary Pieri's rule for Schur polynomials.
So $c^{w'}_{w\, r[k,m]} = 1$.
\QED

Theorem \ref{thm:main}$'$ determines the structure constants
$c^{w'}_{w\, r[k,m]}$ and $c^{w'}_{w\, c[k,m]}$.
We compute more structure constants.
For $\nu$ a partition with $k$ parts, let $w(\nu)$ be the
Grassmannian permutation of descent $k$ and shape $\nu$.

\begin{thm}
Let $w, w'\in S_n$ and $k\leq n$ be an integer.
Suppose  $w\leq_k w'$ and $\ell(w') = \ell(w) +m$.
Let $a_1,b_1,\ldots,a_m,b_m$ be such that
$a_i\leq k <b_i$ where $w' = w t_{a_1\,b_1}\cdots t_{a_m\,b_m}$
and $\ell(w t_{a_1\,b_1}\cdots t_{a_i\,b_i}) = \ell(w) +i$.
Let $\nu$ be a partition with $k$ parts.
\begin{enumerate}
\item  If $\,w\rkm w'$,
the structure constant $c^{w'}_{w\, w(\nu)}$
equals the Littlewood-Richardson coefficient
$c^\lambda_{\mu\,\nu}$, where $\lambda/\mu$ is a skew row of length
$m$ whose $j$th  row has length
$\#\{i \,|\, a_i = j\}$.
\item
If $\,w\ckm w'$,  the structure constant $c^{w'}_{w\, w(\nu)}$
equals the Littlewood-Richardson coefficient
$c^\lambda_{\mu\,\nu}$, where $\lambda/\mu$ is a skew column of length
$m$ whose $j$th  column has length
$\#\{i \,|\, b_i=j\}$.
\end{enumerate}
\end{thm}

\noindent{\bf Proof:} Using the involution
$\frak{S}_{w} \mapsto \frak{S}_{w_0ww_0}$, it suffices to
prove part (1).
We use  part (2) of Lemma~\ref{lemma:pushforward} to
evaluate $c^{w'}_{w\,w(\nu)}$.
Recall that $\frak{S}_{w(\nu)} = \pi^*(s_\nu)$.  Then
\begin{eqnarray*}
c^{w'}_{w\,w(\nu)}\, s_{(n-k)^k}\   =\
\pi_*(c^{w'}_{w\,w(\nu)}\,  \frak{S}_{w_0}) &=&
\pi_*(\frak{S}_w\cdot\frak{S}_{w_0w'}\cdot \frak{S}_{w(\nu)})\\
&=& \pi_*(\frak{S}_w\cdot\frak{S}_{w_0w'}) \cdot s_\nu\\
&=& s_\mu \cdot s_{\lambda^c}\cdot s_\nu\\
&=& c^{\lambda}_{\mu\nu}\,  s_{(n-k)^k}. \ \ \ \  \QED
\end{eqnarray*}

The formulas of Theorem~\ref{thm:main} may be formulated as the
sum over certain paths in the $k$-Bruhat order.
We explain this formulation here.
A (directed) path in the $k$-Bruhat order from $w$ to $w'$
is equivalent to a choice of integers $a_1,b_1,\ldots, a_m,b_m$
with $a_i\leq k < b_i$ and if $w^{(0)} = w$ and
 $w^{(i)} = w^{(i-1)}\cdot t_{a_i\,b_i}$,
then $\ell(w^{(i)}) = \ell(w) + i$ and $w^{(m)} = w'$.
In this case the path is
$$
w = w^{(0)}<_k w^{(1)} <_k w^{(2)} <_k \cdots <_k w^{(m)} = w'.
$$

\begin{lemma}
Let $w, w' \in S_n$ and $k,m$ be positive integers.
Then
\begin{enumerate}
\item $w\rkm w'$ if and only if there is a path in the $k$-Bruhat
order of length $m$ such that
$$
w^{(1)}(a_1) < w^{(2)}(a_2) < \cdots < w^{(m)}(a_m).
$$
\item $w\ckm w'$ if and only if there is a path in the $k$-Bruhat
order of length $m$ such that
$$
w^{(1)}(a_1) > w^{(2)}(a_2) > \cdots > w^{(m)}(a_m).
$$
\end{enumerate}
Furthermore, these paths are unique.
\end{lemma}

 \noindent{\bf Proof:}
If $w\rkm w'$, one may show that
the set of values $\{ w^{(i)}(a_i)\}$ and the set of transpositions
$\{t_{a_i\,b_i}\}$ depend only upon $w$ and $w'$, and not on the
particular path chosen from $w$ to $w'$ in the $k$-Bruhat order.

It is also the case that
rearranging the set $\{ w^{(i)}(a_i)\}$ in order, as in (1),
may be accomplished by interchanging transpositions $t_{a_i\,b_i}$ and
$t_{a_j\,b_j}$ where $a_i\neq a_j$ (necessarily $b_i\neq b_j$).
Both (1) and the uniqueness of this representation follow from these
observations.
Statement (2) follows for similar reasons.
\QED

For a path $\gamma$ in the $k$-Bruhat order, let
$\mbox{end}(\gamma)$ be the endpoint of $\gamma$.
We state a reformulation of Theorem 1.

\begin{cor}[Path formulation of Theorem 1]\label{cor:pieri_paths}
Let $w\in S_n$.
\begin{enumerate}
\item $ \frak{S}_w \cdot \frak{S}_{r[k,m]} \ =\
\sum_\gamma \frak{S}_{\mbox{\scriptsize end}(\gamma)},
$
the sum over all paths $\gamma$ in the $k$-Bruhat order which
start at $w$ such that
$$
w^{(1)}(a_1) < w^{(2)}(a_2) < \cdots < w^{(m)}(a_m),
$$
where $\gamma$ is the path $w <_k w^{(1)} <_k w^{(2)} <_k \cdots
<_k w^{(m)}$.

Equivalently, $c^{w'}_{w\,r[k,m]}$ counts the number of paths
 $\gamma$ in the $k$-Bruhat order  which start at $w$ such that
$$
w^{(1)}(a_1) < w^{(2)}(a_2) < \cdots < w^{(m)}(a_m).
$$
\item $ \frak{S}_w \cdot \frak{S}_{c[k,m]} \ =\
\sum_\gamma \frak{S}_{\mbox{\scriptsize end}(\gamma)},
$
the sum over all paths $\gamma$ in the $k$-Bruhat order which start at
$w$ such that
$$
w^{(1)}(a_1) > w^{(2)}(a_2) > \cdots > w^{(m)}(a_m),
$$
where $\gamma$ is the path $w<_k w^{(1)} <_k w^{(2)} <_k \cdots
<_k w^{(m)}$.

Equivalently, $c^{w'}_{w\,r[k,m]}$ counts  the number of paths
 $\gamma$ in the $k$-Bruhat order which start at
$w$ such that
$$
w^{(1)}(a_1) > w^{(2)}(a_2) > \cdots > w^{(m)}(a_m).
$$
\end{enumerate}
\end{cor}

This is the form of the conjectures of
Bergeron and Billey~\cite{Bergeron_Billey}, and it
exposes a link between multiplying
Schubert polynomials and paths in the Bruhat order.
Such a link is not unexpected.
The Littlewood-Richardson rule
for multiplying Schur functions may be expressed
as a sum over certain paths in Young's lattice of
partitions.
A connection between paths in the Bruhat order and the
intersection theory of Schubert varieties is described
in~\cite{Hiller_intersections}.
We believe the eventual description of the structure
constants $c^w_{uv}$ will be in terms
of counting paths of certain types in the  Bruhat order
on $S_n$, and will yield new results about the Bruhat order on $S_n$.
Corollary~\ref{cor:hook_enumeration} below is one such result.
\smallskip

Using multiset notation for partitions,  $(p,1^{q-1})$
is the hook shape partition whose Young diagram
is the union of a row of length $p$ and
a column of length $q$.
Define $h[k;\,p,q]$ to be the Grassmannian permutation of
descent $k$  and shape $(p,1^{q-1})$.
Then $\frak{S}_{h[k;\,p,q]} = \pi^* s_{(p,1^{q-1})}$.
This permutation, $h[k;\,p,q]$, is the $p+q$-cycle
$$
(k{-}q{+}1\,\,\,k{-}q{+}2\,\ldots\,k{-}1\,\,\,
k\,\,\,k{+}p\,\,\,k{+}p{-}1\,\ldots\,k{+}1).
$$

\begin{thm}\label{thm:hook_formula}
Let $q\leq k$ and $k{+}p \leq n$ be integers.
Set $m = p{+}q{-}1$.
For $w\in S_n$,
$$
\frak{S}_w \cdot \frak{S}_{h[k;\,p,q]} \   =\
\sum \frak{S}_{end(\gamma)},
$$
the sum over all paths $\gamma: w <_k w^{(1)} <_k w^{(2)} <_k \cdots
<_k w^{(m)}$ in the $k$-Bruhat order  with
$$
w^{(1)}(a_1) < \cdots < w^{(p)}(a_p)
\ \ \ \mbox{and}\ \ \
w^{(p)}(a_p) > w^{(p{+}1)}(a_{p{+}1})>\cdots > w^{(m)}(a_m).
$$
Alternatively, those paths $\gamma$  with
$$w^{(1)}(a_1) > \cdots > w^{(q)}(a_q)
\ \ \ \mbox{and}\ \ \
w^{(q)}(a_q) <\cdots < w^{(m)}(a_m).
$$
\end{thm}

Setting either $p=1$ or $q=1$, we recover Theorem~\ref{thm:main}.
If we consider the coefficient $c^{w'}_{w\,h[k;p,q]}$
 of $\frak{S}_{w'}$ in the
product $\frak{S}_w \cdot \frak{S}_{h[k;p,q]}$, we obtain:

\begin{cor}\label{cor:hook_enumeration}
Let $w, w' \in S_n$, and $p,q$ be positive integers where
$\ell(w')-\ell(w) = p+q-1 = m$.
Then the number of paths
$w <_k w^{(1)} <_k w^{(2)} <_k \cdots
<_k w^{(m)} = w'$
in the $k$-Bruhat order from $w$ to $w'$ with
$$
w^{(1)}(a_1) < \cdots < w^{(p)}(a_p)
\ \ \ \mbox{and}\ \ \
w^{(p)}(a_p) > w^{(p{+}1)}(a_{p{+}1})>\cdots > w^{(m)}(a_m)
$$
equals the number of paths with
$$w^{(1)}(a_1) > \cdots > w^{(q)}(a_q)
\ \ \ \mbox{and}\ \ \
w^{(q)}(a_q) <\cdots < w^{(m)}(a_m).
$$
\end{cor}

\noindent{\bf Proof of Theorem~\ref{thm:hook_formula}:}
By the classical Pieri's rule,
$$
s_{p}\, \cdot \, s_{1^{(q-1)}} \   = \
s_{(p{+}1,1^{q-2})} \, + \, s_{(p,1^{q-1})}.
$$
Expressing these as Schubert polynomials (applying $\pi^*$), we have:
$$
\frak{S}_{r[k,p]}\, \cdot\, \frak{S}_{c[k,q{-}1]}
\   =\   \frak{S}_{h[k;\,p{+}1,q{-}1]}  +
\frak{S}_{h[k;\,p,q]}.
$$

Induction on either
$p$ or $q$ (with $m$ fixed) and
Corollary~\ref{cor:pieri_paths} completes the proof.
\QED

\section{Geometry of Intersections}

We deduce Lemma~\ref{lemma:pushforward} by studying certain
intersections of Schubert varieties.
A key fact we use is that if
$X_w\Fdot$ and $X_v\Gdot$ intersect generically transversally,
then
$$
[X_w\Fdot\bigcap X_v\Gdot] \ = \ [X_w\Fdot]\cdot[X_v\Gdot]
 \ = \ \frak{S}_w \cdot \frak{S}_v
$$
in the cohomology ring.
Flags $\Fdot$ and $\Gdot$ are {\em opposite} if for
$1\leq i \leq n$, $F_i + G_{n-i} = V$.
The set of pairs of opposite flags form the dense orbit of the
general linear group $GL(V)$ acting on the space of all pairs
of flags.
Using this observation and Kleiman's Theorem concerning the
transversality of a general translate~\cite{Kleiman},
we conclude that for any $w,v\in S_n$ and opposite flags $\Fdot$ and
$\Gdot$,
$X_w\Fdot$ and $X_v\Gdot$ intersect generically transversally.
(One may also check this directly by examining the tangent spaces.)
In this case the intersection is either empty or it is
irreducible and  contains a
dense subset isomorphic to $(\Bbb{C}^\times\!)^m$,
where $m+ \ell(w) + \ell(v)= \frac{1}{2}n(n-1)$ (cf.~\cite{Deodhar}).
These facts hold for the Schubert subvarieties of $G_kV$ as well.
Namely, if  $\lambda$ and $\mu$ are any partitions and
$\Fdot$ and $\Gdot$ are opposite flags,
then $\Omega_\lambda\Fdot \bigcap \Omega_\mu \Gdot$ is either
empty or it is an irreducible,
generically transverse intersection containing a dense subset isomorphic
to $(\Bbb{C}^\times\!)^m$,
where $m+ |\lambda|+|\mu| = k(n-k)$.

Let $\Fdot$ and $\Fpdot$ be opposite flags in $V$.
Let $e_1,\ldots,e_n$ be a basis for $V$ such that
$e_i$ generates the one dimensional subspace
$F_{n+1-i}\bigcap F'_i$.
We deduce Lemma~\ref{lemma:pushforward} from the following two
results of this section.

\begin{lemma}\label{lemma:geometry_statementI}
Let $w, w' \in S_n$ with $w <_k w'$ and
$\ell(w')-\ell(w) =m$.
Suppose that
$w' = w t_{a_1\,b_1}\cdots t_{a_m\,b_m}$ with
$a_i\leq k <b_i$ for $1\leq 1\leq m$
and $\ell(w t_{a_1\,b_1}\cdots t_{a_i\,b_i}) = \ell(w)+i$.
Let $\pi : \Bbb{F}(V) \rightarrow G_kV$ be the canonical projection.
Define $Y=\langle e_{w(j)}\,|\,j\leq k\,\mbox{ or }\,w(j)\neq w'(j)\rangle$.
Then $Y$ has codimension $d=n- k - \#\{b_1,\ldots,b_m\}$
and
$$
\pi ( X_w \Fdot \bigcap X_{w_0w'}\Fpdot ) \subset G_k Y.
$$
Also, if $\Edot=\Span{f_1,\ldots,f_n}
\in X_w \Fdot \bigcap X_{w_0w'}\Fpdot$, then we may assume that for
$j>k$ with $w(j) = w'(j)$, we have $f_j = e_{w(j)}$.

\end{lemma}

\begin{lemma}\label{lemma:geometry_statementII}
Let $w, w' \in S_n$
with $w \rkm w'$ and let
$a_1,\ldots,b_m$ be as in the statement of
Lemma~\ref{lemma:geometry_statementI}.
Then there exist opposite flags
$\Gdot$ and $\Gpdot$ and partitions
$\lambda\supset\mu$, with $\lambda/\mu$ a skew row of length $m$
whose $j$th row has  length $\#\{i\,|\, a_i = j\}$
such that
$$
\pi ( X_w \Fdot \bigcap X_{w_0w'}\Fpdot )
\quad =\quad \Omega_{\mu}\Gdot \bigcap \Omega_{\lambda^c}\Gpdot,
$$
and the map $\pi|_{X_w \Fdot \bigcap X_{w_0w'}\Fpdot}
:  X_w \Fdot \bigcap X_{w_0w'}\Fpdot \rightarrow
\Omega_{\mu}\Gdot \bigcap \Omega_{\lambda^c}\Gpdot$
has degree 1.
\end{lemma}

Lemma~\ref{lemma:geometry_statementII} is the surprising connection
to the classical Pieri's rule that was mentioned in the Introduction.
A typical geometric proof of Pieri's rule for Grassmannians
(see~\cite{Griffiths_Harris,Hodge_Pedoe}) involves showing a
triple intersection of Schubert varieties
\begin{equation}\label{eq:triple_intersection}
\Omega_{\lambda}\Gdot \bigcap \Omega_{\mu^c}\Gpdot
\,\bigcap\, \Omega_m\Gppdot
\end{equation}
is transverse and consists of a single point,
when $\Gdot, \Gpdot$, and $\Gppdot$ are in suitably general position.

We would like to construct a proof of Theorem~\ref{thm:main} along those lines,
studying a triple intersection of Schubert subvarieties
\begin{equation}\label{eq:triple_intersection_II}
X_w \Gdot \bigcap X_{w_0w'}\Gpdot
\,\bigcap\, X_{r[k,m]}\Gppdot,
\end{equation}
where $\Gdot, \Gpdot$, and $\Gppdot$ are in  suitably general position.
Doing so, one observes that the geometry of the intersection
of~(\ref{eq:triple_intersection_II}) is governed entirely by the
geometry of an intersection similar to that in~(\ref{eq:triple_intersection}).
In part, that is because $ X_{r[k,m]}\Gppdot=\pi^{-1}\Omega_m\Gppdot$.
This is the spirit of our method, which may be seen most vividly in
Lemmas 14 and 15.
\bigskip

\noindent{\bf Proof of Lemma~\ref{lemma:pushforward}:}
Since $\Fdot$ and $\Fpdot$ are opposite flags,
$X_w \Fdot \bigcap X_{w_0w'}\Fpdot$
is a generically transverse intersection, so in the cohomology ring
$$
[X_w \Fdot \bigcap X_{w_0w'}\Fpdot]\  = \
[X_w \Fdot]\cdot [X_{w_0w'}\Fpdot]\  =\  \frak{S}_w\cdot\frak{S}_{w_0w'}.
$$
Let $Y$ be the subspace of Lemma~\ref{lemma:geometry_statementI}.
Since $\pi ( X_w \Fdot \bigcap X_{w_0w'}\Fpdot ) \subset G_k Y$,
the class $\pi_*(\frak{S}_w\cdot\frak{S}_{w_0w'})$ is a cohomology class
on $G_kY$.
However, all such classes are of the form $\delta \cdot [ G_k Y]$,
for some cohomology class $\delta$ on $G_kV$.
Since $d$ is the codimension of $Y$, we have $[G_k Y]= s_{d^k}$,
establishing part (1) of Lemma~\ref{lemma:pushforward}.

For part (2), suppose further that $w \rkm w'$.
If $\rho$ is the restriction of $\pi$ to
$X_w \Fdot \bigcap X_{w_0w'}\Fpdot$, then
$$
\pi_*(\frak{S}_w\cdot\frak{S}_{w_0w'}) \ = \
\pi_*([X_w \Fdot \bigcap X_{w_0w'}\Fpdot])
 \ = \ \deg \rho\cdot [\pi(X_w \Fdot \bigcap X_{w_0w'}\Fpdot)].
$$
By Lemma~\ref{lemma:geometry_statementII},
$\deg \rho = 1$ and
$\pi(X_w\Fdot \bigcap x_{w_0w'}\Fpdot) =
\Omega_\mu \Gdot \bigcap \Omega_{\lambda^c} \Gpdot$.
Since
$\Gdot$ and $\Gpdot$ are opposite flags, we have
$$
\pi_*(\frak{S}_w\cdot\frak{S}_{w_0w'}) \ = \
1\cdot [ \Omega_{\lambda}\Gdot \bigcap \Omega_{\mu^c}\Gpdot]
\ = \ [ \Omega_{\lambda}\Gdot ] \cdot [\Omega_{\mu^c}\Gpdot]
\ =\  s_\lambda \cdot s_{\mu^c},
$$
completing the proof of Lemma~\ref{lemma:pushforward}.\QED
\bigskip

We deduce Lemma~\ref{lemma:geometry_statementI} from a series of lemmas.
We first make a definition.
Let $W\subsetneq V$ be a codimension 1 subspace and let $e \in V - W$
so that $V = \Span{W,e}$.
For $1\leq p \leq n$, define an expanding map
$\psi_p: \Bbb{F}(W) \rightarrow \Bbb{F}(V)$ as follows
$$
(\psi_p \Edot)_i = \left\{
\begin{array}{ll} E_i & \mbox{ if } i<p\\
\Span{E_{i-1},e} & \mbox{ if } i \geq p
\end{array}\right.  .
$$
Note that if $\Edot = \Span{f_1,\ldots,f_{n-1}}$, then
$\psi_p \Edot= \Span{f_1,\ldots,f_{p-1},e,f_p,\ldots,f_{n-1}}$.

For $w\in S_n$ and $1\leq p \leq n$, define $w|_p \in S_{n-1}$ by
$$
w|_p (j) = \left\{
\begin{array}{ll}
w(j)       & \mbox{ if } j<p \mbox{ and } w(j)<w(p)\\
w(j{+}1)   & \mbox{ if } j\geq p \mbox{ and } w(j)<w(p)\\
w(j) - 1   & \mbox{ if } j<p \mbox{ and } w(j)>w(p)\\
w(j{+}1)-1 & \mbox{ if } j\geq p \mbox{ and } w(j)>w(p)
\end{array} \right. .
$$
If we represent permutations as matrices, $w|_p$ is obtained
by crossing out the $p$th row and $w(p)$th column of the matrix
for $w$.

\begin{lemma}\label{lemma:expand}
Let $W\subsetneq V$ and $e\in V - W$ with $V = \Span{W,e}$.
Let $\Gdot$ be a complete flag in $W$.
For $1\leq p \leq n$ and
$w\in S_n$,
$$
\psi_p\left( X_{w|_p}\Gdot \right) \subset
X_w \left(\psi_{w_0w(p)}(\Gdot)\right).
$$
\end{lemma}

\noindent{\bf Proof:}
Let $\Edot \in X_{w|_p}\Gdot $.   Then $W$ has a basis
$f_1,\ldots,f_{n-1}$  with $\Edot = \Span{f_1,\ldots,f_{n-1}}$
and for each $1\leq i \leq n-1$,
$f_i \in G_{n-w|_p(i)}$.
Then we necessarily have
$
\psi_p(\Edot) = \Span{\phi_1,\ldots,\phi_n}
= \Span{f_1,\ldots,f_{p-1},e,f_p,\ldots,f_{n-1}}$.
Noting
$$
\left(\psi_{w_0w(p)}(\Gdot)\right)_{n+1-j} \ = \ \left\{
\begin{array}{ll} G_{n+1-j} & \mbox{ if } j >w(p) \\
\Span{e,G_{n-j}}  & \mbox{ if } j \leq w(p)   \end{array}\right. ,
$$
we see that $\phi_i \in
\left(\psi_{w_0w(p)}(\Gdot)\right)_{n+1-w(i)}$.
Thus $\psi_p\left( X_{w|_p}\Gdot \right) \subset
X_w \left(\psi_{w_0w(p)}(\Gdot)\right)$.
\QED

\begin{lemma}\label{lemma:expanding}
Let $W\subsetneq V$ and $e\in V - W$ with $V = \Span{W,e}$ and
let $\Gdot$ and $\Gpdot$ be opposite flags in $W$.
Suppose that $w<_kw'$ are permutations in $S_n$ and
$p>k$ an integer such that $w(p) = w'(p)$.
Let $w_0^{(j)}$ is the longest permutation in $S_j$.
Then
\begin{enumerate}
\item $\ell(w'|_p)-\ell(w|_p)=\ell(w')-\ell(w)$ and $w|_p <_k w'|_p$.
\item $\psi_p{\left(X_{w|_p}\Gdot \bigcap X_{w_0^{(n-1)}(w'|_p)}\Gpdot\right)}
\, = \,X_w\left(\psi_{w_0^{(n)}w(p)}(\Gdot) \right)
\bigcap X_{w_0^{(n)}w'}\left(\psi_{w'(p)}(\Gpdot)\right)$.
\item If $\Edot \in
 X_w\left(\psi_{w_0^{(n)}w(p)}(\Gdot) \right)
\bigcap X_{w_0^{(n)}w'}\left(\psi_{w'(p)}(\Gpdot)\right)$,
then $E_p = \Span{E_{p-1},e}$.
\item If $\Fdot$ and $\Fpdot$ are opposite flags
in $V$  and
$\Edot \in X_w\Fdot\bigcap X_{w_0^{(n)}w'}\Fpdot$,
then $E_k \in F_{n-w(p)}+F'_{w(p)-1}$.
\end{enumerate}
\end{lemma}

\noindent{\bf Proof:}
First recall that $\ell(v t_{a\,b}) = \ell(v)+1$ if and only if
$v(a) < v(b)$ and if $a<j<b$, then
$v(j)$ is not between $v(a)$ and $v(b)$.
Thus if $\ell(v t_{a\,b}) = \ell(v)+1$ and $p\not\in \{a,b\}$,
we have  $\ell(v t_{a\,b}|_p) = \ell(v|_p)+1$.
Statement (1) follows by induction on
$\ell(w') - \ell(w)$.

For (2), since $(w_0^{(n)}w')|_p = w_0^{(n-1)}(w'|_p)$ and
$w_0^{(n)}w_0^{(n)}w'= w'$, Lemma~\ref{lemma:expand}
shows
$$
\psi_p\left(X_{w|_p}\Gdot \bigcap
X_{w_0^{(n-1)}(w'|_p)}\Gpdot\right)
\subset X_w\left(\psi_{w_0^{(n)}w(p)}(\Gdot)\right)
\bigcap X_{w_0^{(n)}w'}\left(\psi_{w'(p)}(\Gpdot)\right).
$$
The flags $\psi_{w_0^{(n)}w(p)}(\Gdot)$ and
$\psi_{w'(p)}(\Gpdot)$ are opposite flags in $V$, since
$\Gdot$ and $\Gpdot$ are opposite flags in $W$.
Then part (1) shows both sides have the same dimension.
Since $\psi_p$ is injective, they are equal.

To show (3), let $\Edot \in
 X_w\left(\psi_{w_0^{(n)}w(p)}(\Gdot) \right)
\bigcap X_{w_0^{(n)}w'}\left(\psi_{w'(p)}(\Gpdot)\right)$.
By (2), there is a flag
$\Epdot \in X_{w|_p}\Gdot \bigcap
X_{w_0^{(n-1)}(w'|_p)}\Gpdot$
with $\psi_p(\Epdot) = \Edot$,
so $E_p = \Span{E'_{p-1},e} =  \Span{E_{p-1},e}$.

For (4), let $W= F_{n-w(p)}+F'_{w'(p)-1}$ and $e$
any nonzero vector in the one dimensional space
$F_{n+1-w(p)}\bigcap F'_{w'(p)}$.
The distinct subspaces in $\Fdot\bigcap W$ define a flag $\Gdot$,
and those in $\Fpdot \bigcap W$ define a flag $\Gpdot$.
In fact,  $\psi_{w_0^{(n)}w(p)}(\Gdot)=\Fdot$
and $\psi_{w(p)}(\Gpdot) = \Fpdot$, and $\Gdot$ and
$\Gpdot$ are opposite flags in $W$.
By (2),
$$
\psi_p\left(X_{w|_p}\Gdot \bigcap
X_{w_0^{(n-1)}(w'|_p)}\Gpdot\right)
= X_w\Fdot
\bigcap X_{w_0^{(n)}w'}\Fpdot.
$$
Thus flags in
$X_w\Fdot \bigcap X_{w_0^{(n)}w'}\Fpdot$
are in the image of $\psi_p$.
As $k < p$,  $\left( \psi_p \Edot\right)_k = E_k\subset W$,
establishing part (4).
\QED

\noindent{\bf Proof of Lemma~\ref{lemma:geometry_statementI}:}
Let $\Fdot$ and $\Fpdot$ be opposite flags in $V$, let
$w<_k w'$ and let
$\Edot \in X_w\Fdot \bigcap X_{w_0w'}\Fpdot$.
Define a basis $e_1,\ldots,e_n$ for $V$ by
$F_{n+1-j}\bigcap F'_j = \Span{e_j}$ for $1\leq j\leq n$.
Suppose $w' = w t_{a_1\,b_1}\cdots t_{a_m\,b_m}$ with
$a_i\leq k < b_i$.
Let $\{p_1,\ldots,p_d\}$ be the complement of
$\{b_1,\ldots,b_m\}$ in $\{k+1,\ldots,n\}$.
For $1\leq i\leq d$, let
$Y_i = \Span{e_1,\ldots,\widehat{e_{w(p_i)}},\ldots,e_n}
 = \Span{e_1,\ldots,e_{w(p_i)-1},e_{w(p_i)+1},\ldots,e_n}$.
Since $w(p_i) = w'(p_i)$ and $k<p_i$, we see that
$Y_i = F_{n-w(p_i)} + F'_{w(p_i)-1}$, so part (4) of
Lemma~\ref{lemma:expanding} shows $E_k \subset Y_i$.
Thus
$$
E_k \in
\bigcap_{i=1}^d Y_i \ = \ \Span{e_{w(j)}\,|\, j<k\mbox{ or }j=b_i}
\  = \ Y.
$$
Since $w(p_i) = w'(p_i)$ for $1\leq i \leq d$,
we have $E_{p_i} = \Span{E_{p_i-1},e_{w(p_i)}}$, by  part (3) of
Lemma~\ref{lemma:expanding}.
So if $\Edot = \Span{f_1,\ldots,f_n}$, we may assume that
$f_{p_i} = e_{w(p_i)} \in F_{n+1-w(p_i)} \cap F'_{w'(p_i)}$ for
$1\leq i \leq d$, completing the proof.
\QED

To prove Lemma~\ref{lemma:geometry_statementII}, we begin by describing
an intersection in a  Grassmannian.
Recall that $\Omega_\lambda\Fdot =
\{H\in G_kV\,|\, \dim H\cap F_{k-j+\lambda_j} \geq j \mbox{ for }
1\leq j\leq k\}$.

\begin{lemma}\label{lemma:grassmannian}
Suppose that $L_1,\ldots,L_k,M \subset V$ with
$V = M  \bigoplus L_1\bigoplus\cdots\bigoplus L_k$.
Let $r_j= \dim L_j -1$ and $m = r_1 + \cdots + r_k$.
Then there are opposite flags $\Fdot$ and $\Fpdot$ and
partitions $\lambda\supset \mu$ with $\lambda_j - \mu_j = r_j$
and $\lambda/\mu$ a skew row of length $m$ such that in
$G_kV$,
$$
\Omega_\mu\Fdot \bigcap \Omega_{\lambda^c}\Fpdot =
\{ H\in G_kV\,|\, \dim H\bigcap L_j = 1 \mbox{ for } 1\leq j \leq k\}.
$$
\end{lemma}

\noindent{\bf Proof:}
Let $\mu_k =0$ and $\mu_j = r_{j+1}+\cdots+r_k$ for
$1\leq j <k$ and $\lambda_j =r_j+\mu_j$ for $1\leq j\leq k$.
Choose a basis $e_1,\ldots,e_n$ for $V$ such that
\begin{eqnarray*}
L_j & =&
\Span{e_{k+1-j+\mu_j},e_{k+2-j+\mu_j},\ldots,e_{k+1+r_j-j+\mu_j}
=e_{k+1-j+\lambda_j}}\\
M &=& \Span{e_{m+k+1},\ldots,e_n}
\end{eqnarray*}
Let $\Fdot = \Span{e_n\ldots,e_1}$ and $\Fpdot=\Span{e_1,\ldots,e_n}$.
Then
\begin{center}
$F_{n-k+j-\mu_j}\quad  =\quad
M\bigoplus L_1 \bigoplus \cdots \bigoplus L_j \quad\quad$\\
$F'_{n-k+(k+1-j)-\lambda^c_{k+1-j}} \quad = \quad F'_{k+1-j+\lambda_j}
\quad  =\quad
L_j\bigoplus \cdots \bigoplus L_k$.
\end{center}
If $H\in \Omega_\mu\Fdot \bigcap \Omega_{\lambda^c}\Fpdot$, then
$\dim H\bigcap F_{n-k+j-\mu_j} \geq j$ for $1\leq j\leq k$
and
$$
\dim H\bigcap F'_{n-k+(k+1-j)-\lambda^c_{k+1-j}}\geq k+1-j,
$$
for $1\leq j\leq k$.
Thus for $1\leq j \leq k$,
$$
\dim H\bigcap F_{n-k+j-\mu_j} \bigcap
F'_{n-k+(k+1-j)-\lambda^c_{k+1-j}} \geq 1.
$$
But $F_{n-k+j-\mu_j} \bigcap F'_{n-k+(k+1-j)-\lambda^c_{k+1-j}} = L_j$,
so $\dim H\bigcap L_j\geq 1$ for $1\leq j \leq k$.
Since $L_j\bigcap L_i = \{0\}$ if $j\neq i$, we see that
$\dim H\bigcap L_j =1$.
Thus
$$
\Omega_\mu\Fdot \bigcap \Omega_{\lambda^c}\Fpdot \subset
\{ H\in G_kV\,|\, \dim H\bigcap L_j = 1 \mbox{ for } 1\leq j \leq k\}.
$$
We show these varieties have the same dimension,
establishing their equality:
Since $|\lambda| = |\mu| +m$, and $\Fdot$ and $\Fpdot$ are opposite
flags, $\Omega_\mu\Fdot \bigcap \Omega_{\lambda^c}\Fpdot$
has dimension $m$.
But  the map
$H \mapsto (H\bigcap L_1,\ldots,H\bigcap L_k)$
defines an isomorphism between
$\{ H\in G_kV\,|\, \dim H\bigcap L_j = 1 \mbox{ for } 1\leq j \leq k\}$
and $\Bbb{P}L_1 \times \cdots \times \Bbb{P}L_k$,
which  has dimension $\sum_j (\dim L_j-1) = m$.
Here, $\Bbb{P}L_j$ is the projective space of one dimensional
subspaces of $L_j$.
\QED

We relate this to intersections of Schubert varieties in the
flag manifold.

\begin{lemma}\label{lemma:intersection_calculation}
Suppose that $w\rkm w'$ and
$w' = wt_{a_1\,b_1}\cdots t_{a_m\,b_m}$ with $a_i\leq k<b_i$ and
$\ell(wt_{a_1\,b_1}\cdots t_{a_i\,b_i}) = \ell(w)+i$.
Let $\Fdot$ and $\Fpdot$ be opposite flags in $V$ and let
$\Span{e_i} = F_{n+1-i}\bigcap F'_{i}$.
Define
\begin{eqnarray*}
L_j & =& \Span{e_j, e_{w(b_i)}\,|\, a_i = j}\\
M  & =&\Span{e_{w(p)}\,|\, k<p \mbox{ and } w(p) = w'(p)}.
\end{eqnarray*}
Then
\begin{enumerate}
\item  $\dim L_j = 1 +\#\{i\,|\, a_i = j\}$ and
$ V = M \bigoplus L_1\bigoplus \cdots \bigoplus L_k$.
\item If $\Edot \in X_w\Fdot \bigcap X_{w_0w'}\Fpdot$, then
$\dim E_k \bigcap L_j = 1$ for $1\leq j\leq k$.
\item
Let $\pi$ be the map induced by $\Edot \mapsto E_k$.
Then
$$
\pi : X_w\Fdot \bigcap X_{w_0w'}\Fpdot
\rightarrow \{H\in G_kV\,|\, \dim H\bigcap L_j = 1 \mbox{ for }
1\leq j\leq k\}
$$
is surjective and of degree 1.
\end{enumerate}
\end{lemma}

\noindent{\bf Proof:}
Part (1) is immediate.

For (2) and (3), note that both
$\{H\in G_kV\,|\, \dim H\bigcap L_j = 1 \mbox{ for }
1\leq j\leq k\}$ and $X_w\Fdot \bigcap X_{w_0w'}\Fpdot$
are irreducible and have dimension $m$.
We exhibit an $m$ dimensional subset of each over which $\pi$ is
an isomorphism.

Let $\alpha = (\alpha_1,\ldots,\alpha_m) \in (\Bbb{C}^\times\!)^m$
be an $m$-tuple of nonzero complex numbers.
We define a basis $f_1,\ldots,f_n$ of $V$ depending upon
$\alpha$ as follows.
$$
f_j = \left\{
\begin{array}{ll}
e_{w(j)} +
{\displaystyle \sum_{i : a_i=j}
\alpha_i e_{w(b_i)}} &
\mbox{ if } j\leq k \\
e_{w(j)} & \mbox{ if } j>k \mbox{ and } j\not\in\{b_1,\ldots,b_m\}\\
\rule{0pt}{22pt}{\displaystyle
\sum_{\shortstack{\scriptsize $i: a_i = a_q$\\ \scriptsize $w(b_i)\geq w(j)$}}
\alpha_i e_{w(b_i)}} & \mbox{ if } j = b_q > k
\end{array} \right. .
$$
Let $i_1<\cdots < i_s$ be those
integers $i_l$ with $a_{i_l} = j$.
Since $t_{a_i\,b_i}$ lengthens the permutation
$wt_{a_1\,b_1}\cdots t_{a_{i-1}\,b_{i-1}}$, we see that
$$
\begin{array}{ccccccc}
w(j)         & < &  w(b_{i_1}) & < & \cdots & < &  w(b_{i_s}) \\
\parallel    &   &  \parallel  &   &        &   &  \parallel  \\
w'(b_{i_1})  & < &  w'(b_{i_2}) & < & \cdots & < &     w'(j)
\end{array}
$$
Thus the first term in $f_j$ is proportional to $e_{w(j)}$.
Hence $f_j \in F_{n+1-w(j)} - F_{n-w(j)}$, and so
$f_1,\ldots,f_n$ is a basis  of $V$ and the flag
$\Edot(\alpha) = \Span{f_1,\ldots,f_n}$ is in $X_w\Fdot$.

Note that $f'_1,\ldots,f'_n$ is also a basis for $\Edot(\alpha)$, where
$f'_j$ is given by
$$
f'_j = \left\{
\begin{array}{ll}
f_j &
\mbox{ if } j\leq k \\
f_j& \mbox{ if } j>k \mbox{ and } j\not\in\{b_1,\ldots,b_m\}\\
f_{a_q} - f_j & \mbox{ if } j = b_q > k
\end{array} \right. .
$$
Here, the last term in each $f'_j$ is proportional to $e_{w'(j)}$,
so $f'_j \in F'_{w'(j)} = F'_{n+1-w_0w'(j)}$,
showing that
$\Edot(\alpha) \in X_{w_0w'}\Fpdot$.

Since  $f_j \in L_j$ for $1\leq j\leq k$, we have
$\dim \Edot(\alpha) \bigcap L_j = 1$ for $1\leq j\leq k$.
As $\{\Edot(\alpha)\,|\, \alpha \in (\Bbb{C}^\times\!)^m\}$
is a subset of  $ X_w\Fdot\bigcap X_{w_0w'}\Fpdot$ of dimension $m$,
it is dense.
Thus if $\Edot \in X_w\Fdot\bigcap X_{w_0w'}\Fpdot$, then
$\dim E_k \bigcap L_j = 1$  for $1\leq j\leq k$.

The set $\{(\Edot(\alpha))_k\,|\, \alpha \in (\Bbb{C}^\times\!)^m\}$
is a dense subset of
$$
\{H\in G_kV\,|\, \dim H\bigcap L_j = 1 \mbox{ for } 1\leq j\leq k\}
\ \simeq \ \Bbb{P}L_1\times\cdots\times\Bbb{P}L_k.
$$
Since $\pi$ is an isomorphism of this set with
$\{\Edot(\alpha)\,|\, \alpha \in (\Bbb{C}^\times\!)^m\}$,
the map
$$
\pi:  X_w\Fdot\bigcap X_{w_0w'}\Fpdot \rightarrow
\{H\in G_kV\,|\, \dim H\bigcap L_j = 1 \mbox{ for } 1\leq j\leq k\}
$$
is surjective of degree 1, proving the lemma.
\QED

We note that Lemma~\ref{lemma:geometry_statementII} is an
immediate consequence of Lemmas~\ref{lemma:grassmannian}
and~\ref{lemma:intersection_calculation}(3).

\section{Examples}

In this section we describe two examples, which should serve to illustrate
the results of Section 5.
This manuscript differs from the version we are submitting
for publication only by the inclusion of this section, and
its mention  in the Introduction.
\bigskip

Fix a basis $e_1,\ldots,e_7$ for $\Bbb{C}\,^7$.
This gives coordinates for vectors in $\Bbb{C}\,^7$, where
$(v_1,\ldots,v_7)$ corresponds to
$v_1e_1{+}\cdots{+}v_7e_7$.
Define the opposite flags
$\Fdot$ and $\Fpdot$ by
$$
\Fdot = \langle e_7,e_6,e_5,e_4,e_3,e_2,e_1\rangle \  \
\mbox{and} \ \
\Fpdot = \langle e_1,e_2,e_3,e_4,e_5,e_6,e_7\rangle.
$$
For example, $F_3 = \langle e_7,e_6,e_5 \rangle$ and
$F'_4 =  \langle e_1,e_2,e_3,e_4 \rangle$.
Let $w= 5412763$, $w' = 6524713$ and $w'' = 7431652$ be permutations in
$S_7$.
(We denote permutations by the sequence of their values.)
Their lengths are 10, 14, and 14, respectively, and
$w<_4 w'$ and $w<_3 w''$.
We seek to describe the intersections
$$
X_w\Fdot \bigcap X_{w_0w'}\Fpdot \ \ \ \
\mbox{and} \ \ \ \ X_w\Fdot \bigcap X_{w_0w''}\Fpdot.
$$
Rather than describe each in full, we describe a dense subset of
each which is isomorphic to the torus, $(\Bbb{C}^\times\!)^4$.
This suffices for our purposes.

Recall that the Schubert cell $X^\circ_w\Fdot$ is defined to be
$$
X^{\circ}_w \Fdot = \{ \Edot=\Span{f_1,\ldots,f_7}\,|\,
f_i \in F_{8-w(i)}-F_{7-w(i)}, \,1\leq i\leq 7\}.
$$
Using the given coordinates of $\Bbb{C}\,^7$,
we may write a typical element of
$X^\circ_w\Fdot$ in a unique manner.
For each $f_i \in F_{8-w(i)}-F_{7-w(i)}$, the coordinate 7-tuple
for $f_i$ has zeroes in the places $1,\ldots,w(i)-1$
and a nonzero coordinate in its $w(i)$th place, which we
assume to be 1.
We may also assume that the $w(j)$th coordinate of $f_i$ is zero for those
$j<i$ with $w(j)> w(i)$, by subtracting a suitable multiple of $f_j$.
Writing the coordinates of $f_1,\ldots,f_7$ as rows of an array, we
conclude that a typical flag in $X^{\circ}_w \Fdot$ has a unique
representation of the following form:
$$
\begin{array}{ccccccc}
\cdot &\cdot&\cdot&\cdot& 1 & * & * \\
\cdot &\cdot&\cdot& 1 &\cdot& * & * \\
1 & * & * &\cdot&\cdot& * & * \\
\cdot & 1 & * &\cdot&\cdot& * & * \\
\cdot &\cdot&\cdot&\cdot&\cdot&\cdot& 1 \\
\cdot &\cdot&\cdot&\cdot&\cdot& 1 &\cdot\\
\cdot &\cdot& 1 &\cdot&\cdot&\cdot& \cdot
\end{array}
$$
Here, the $i$th column contains the coefficients of $e_i$,
the $\cdot$'s represent 0,
and the $*$'s indicate some complex numbers, uniquely
determined by the flag.
Likewise, flags in $X^{\circ}_{w_0w'} \Fpdot$
and $X^{\circ}_{w_0w''} \Fpdot$  have unique bases of the forms:
$$
\begin{array}{cccccccc}
 &  * & * & * & * & * & 1 &\cdot\\
 & * & * & * & * & 1 &\cdot&\cdot\\
 &  * & 1 &\cdot&\cdot&\cdot&\cdot&\cdot\\
 &  * &\cdot& * & 1 &\cdot&\cdot&\cdot\\
 &  * &\cdot& * &\cdot&\cdot&\cdot& 1 \\
 &  1 &\cdot&\cdot&\cdot&\cdot&\cdot&\cdot\\
 & \cdot&\cdot& 1 &\cdot&\cdot&\cdot& \cdot
\end{array}
\hspace{1in}
\begin{array}{cccccccc}
 &  * & * & * & * & * & * & 1 \\
 &  * & * & * & 1 &\cdot&\cdot&\cdot\\
 &  * & * & 1 &\cdot&\cdot&\cdot&\cdot\\
 &  1 &\cdot&\cdot&\cdot&\cdot&\cdot&\cdot\\
 & \cdot& * &\cdot&\cdot& * & 1 &\cdot\\
 & \cdot& * &\cdot&\cdot& 1 &\cdot&\cdot\\
 & \cdot& 1 &\cdot&\cdot&\cdot&\cdot& \cdot
\end{array}
$$

Let $\alpha,\beta,\gamma$ and $\delta$ be four nonzero complex numbers.
Define bases $f_1,f_2,\ldots,f_7$ and $g_1,g_2,\ldots,g_7$ by the
following arrays of coordinates.
$$
\begin{array}{ccccccccc}
f_1 &= &\cdot&  \cdot   &\cdot&  \cdot   &    1  & \alpha &\cdot\\
f_2 &= &\cdot&  \cdot   &\cdot&   1    & \beta &   \cdot  &\cdot\\
f_3 &= & 1 & \gamma &\cdot&  \cdot   &   \cdot &   \cdot  &\cdot\\
f_4 &= &\cdot&   1    &\cdot& \delta &   \cdot &   \cdot  &\cdot\\
f_5 &= &\cdot&  \cdot   &\cdot&  \cdot   &   \cdot &   \cdot  & 1 \\
f_6 &= &\cdot&  \cdot   &\cdot&  \cdot   &   \cdot & \alpha &\cdot\\
f_7 &= &\cdot&  \cdot   & 1 &  \cdot   &   \cdot &   \cdot  &\cdot
\end{array}
\hspace{1in}
\begin{array}{ccccccccc}
g_1 & = &\cdot&  \cdot   &   \cdot  &\cdot& 1 & \alpha & \beta \\
g_2 & = &\cdot&  \cdot   &   \cdot  & 1 &\cdot&   \cdot  &  \cdot  \\
g_3 & = & 1 & \gamma & \delta &\cdot&\cdot&   \cdot  &  \cdot  \\
g_4 & = &\cdot& \gamma & \delta &\cdot&\cdot&   \cdot  &  \cdot  \\
g_5 & = &\cdot&  \cdot   &   \cdot  &\cdot&\cdot&   \cdot  & \beta \\
g_6 & = &\cdot&  \cdot   &   \cdot  &\cdot&\cdot& \alpha & \beta \\
g_7 & = &\cdot&  \cdot   & \delta &\cdot&\cdot&   \cdot  &  \cdot
\end{array}
$$
Let $\Edot = \Span{f_1,f_2,\ldots,f_7}$ and
$\Epdot = \Span{g_1,g_2,\ldots,g_7}$.
Considering the left-most nonzero entry in each row, we see that
both $\Edot$ and $\Epdot$ are in $X^\circ_w\Fdot$.  To see that
$\Edot \in X^\circ_{w_0w'}\Fpdot$ and
$\Epdot \in X^\circ_{w_0w''}\Fpdot$, note that we could
choose
$$
\begin{array}{ccccccccc}
f_6' & = & 1 &  \cdot  &\cdot&\cdot&   \cdot  &   \cdot&\cdot
\end{array}
\hspace{1in}
\begin{array}{ccccccccc}
g_4' &  = & 1 &   \cdot  &\cdot&\cdot&\cdot&   \cdot  &\cdot\\
g_5' &  = &\cdot&   \cdot  &\cdot&\cdot& 1 & \alpha &\cdot\\
g_6' &  = &\cdot&   \cdot  &\cdot&\cdot& 1 &   \cdot  &\cdot\\
g_7' &  = & 1 & \gamma &\cdot&\cdot&\cdot&   \cdot  &\cdot
\end{array}
$$
Replacing the unprimed vectors by the corresponding primed ones
gives alternate bases for $\Edot$ and $\Epdot$.
This shows $\Edot \in X^\circ_{w_0w'}\Fpdot$ and
$\Epdot \in X^\circ_{w_0w''}\Fpdot$.

We use this computation to illustrate Lemmas~\ref{lemma:geometry_statementI}
and~\ref{lemma:geometry_statementII}.

\begin{enumerate}
\item[I.]  First note that for
$\Edot = \Span{f_1,f_2,f_3,f_4,f_5,f_6,f_7}$ as above,
\begin{eqnarray*}
E_3 & \subset & \Span{e_1,e_2,e_5,e_5,e_6}\\
    & = & \Span{ e_{w(j)}\,|\, j\leq k\mbox{ or } w(j)\neq w'(j)}\\
    & = & Y,
\end{eqnarray*}
the subspace of Lemma~\ref{lemma:geometry_statementI}.
Since this holds for all $\Edot$ in a dense subset of
$X_w\Fdot \bigcap X_{w_0w'}\Fpdot$, it holds for all $\Edot$ in
that intersection.

\item[II.]  Recall that $w=5412763$ and note that
$7431652 = w''= w\cdot t_{34}\cdot t_{16}\cdot  t_{37}\cdot t_{15}$,
so $w \stackrel{r[3,4]}{\llra} w''$, and we are in the situation of
Lemma~\ref{lemma:geometry_statementII}.
Let $\mu = (2,2,0)$ and $\lambda = (4,2,2)$ be partitions.
Then $\lambda^c = (2,2,0)$,
and if $\Epdot = \Epdot(\alpha,\beta,\gamma,\delta)$ is a
flag in the above form,
then
$$
E_3'(\alpha,\beta,\gamma,\delta) \in
\Omega_{\mu}\Fdot \bigcap \Omega_{\lambda^c} \Fpdot,
$$
since
$$
\begin{array}{ccccl}
f_1 &\in & F_3 & = & F_{7-3+1-\mu_1} \bigcap F'_{7-3+3-\lambda^c_3}\\
f_2 &\in & \Span{e_4} & =& F_{7-3+1-\mu_2} \bigcap F'_{7-3+3-\lambda^c_2}\\
f_3 &\in & F'_3 & = & F_{7-3+1-\mu_3} \bigcap F'_{7-3+3-\lambda^c_1}.
\end{array}
$$
Furthermore, the map
$\pi: \Epdot \mapsto E_3'$
is injective for those $\Epdot(\alpha,\beta,\gamma,\delta)$ given above.
Since that set is dense in
$X_w\Fdot \bigcap X_{w_0w''}\Fpdot$,
and the set of $E_3'(\alpha,\beta,\gamma,\delta)$ is dense in
$\Omega_{\mu}\Fdot \bigcap \Omega_{\lambda^c} \Fpdot$,
it follows that
$$
\pi : X_w\Fdot \bigcap X_{w_0w''}\Fpdot \rightarrow
\Omega_{\mu}\Fdot \bigcap \Omega_{\lambda^c} \Fpdot
$$
is surjective and of degree 1.
\end{enumerate}
\smallskip

Note that the description of $X_w\Fdot \bigcap X_{w_0w''}\Fpdot$ in II
is consistent with that given for general $w\rkm w''$ in the proof of
Lemma~\ref{lemma:intersection_calculation}, part (2).
This explicit description is the key to the understanding we gained
while trying to establish Theorem~\ref{thm:main}
\smallskip

Also note that $w' = w \cdot t_{16}\cdot t_{26}\cdot t_{46}\cdot t_{36}$,
thus $w \cVIVI w'$.
In I above, we give an explicit description of the intersection
$X_w\Fdot \bigcap X_{w_0w'}\Fpdot$.
This may be generalized to give a similar description whenever
$w \ckm w'$, and may be used to establish Theorem~\ref{thm:main}
in much the same manner as we used the explicit description
of intersections when $w\rkm w'$.

\end{document}